\documentclass{article}

\usepackage[english]{babel}

\usepackage[letterpaper,top=2cm,bottom=2cm,left=3cm,right=3cm,marginparwidth=1.75cm]{geometry}

\usepackage{amsmath}
\usepackage{graphicx}
\usepackage[colorlinks=true, allcolors=blue]{hyperref}
\usepackage{xfrac}       % compact symbols for 1/2, etc

\usepackage{amssymb}
\usepackage{amsthm}
\usepackage{amscd}

\usepackage[version=4]{mhchem} % chemical elements
\usepackage{siunitx}

\usepackage{mathtools}
\usepackage{algorithm}
\usepackage[noend]{algpseudocode}
\usepackage{microtype}

\usepackage{tikz}
\usetikzlibrary{arrows.meta,calc,decorations.markings,fit,positioning,shapes.geometric}

\floatname{algorithm}{Algorithm}

\algrenewcommand\algorithmicrequire{\textbf{Input:}}
\algrenewcommand\algorithmicensure{\textbf{Output:}}

\newcommand{\TTcross}{\operatorname{TT\text{-}Cross}}
\newcommand{\TTround}{\operatorname{TT\text{-}Round}}
\newcommand{\Heat}{\mathcal H_d}

\newcommand{\RelRes}[2]{\dfrac{\lVert #1\rVert_2}
  {\max\{\lVert #2\rVert_2,\varepsilon_{\mathrm{pos}}\}}}

\newcommand{\Id}{\operatorname{I}_d}
\newcommand{\diff}{\mathrm{d}}

\newtheorem{corollary}{Corollary}[section]

\numberwithin{equation}{section}

\title{Bayesian Inverse Problems for EUV reflectometry with regularized JKO scheme}
\author{
    Vitalii Aksenov\thanks{
        Weierstrass Institute for Applied Analysis and Stochastics (WIAS), 
        Address: Anton-Wilhelm-Amo-Str. 39,
        10117 Berlin, Germany.
        Email: \texttt{aksenov@wias-berlin.de}
    }
    \and
    Oliver Henze\thanks{
        Physikalisch-Technische Bundesanstalt (PTB), 
        Address: Abbestraße 2–12, 
        10587 Berlin, Germany.
        Email: \texttt{oliver.henze@ptb.de}
    }
}

\begin{document}
\maketitle

\begin{abstract}
    We study a Bayesian inverse problem in EUV reflectometry, where optical constants and layer thicknesses are inferred from angle- and wavelength-resolved reflection measurements. 
    The physical process is modeled and evaluated using an efficient forward model based on the transfer matrix method. 
    As an alternative to computationally expensive sampling-based inference, we employ an Eulerian, entropically regularized Jordan--Kinderlehrer--Otto (JKO) scheme that minimizes the Kullback--Leibler divergence via proximal steps in the Wasserstein space. 
    Using the Benamou--Brenier dynamic formulation and a Hopf--Cole transform, each JKO step is reformulated as a pair of heat equations coupled only through their initial and terminal conditions, enabling an efficient fixed-point solver. 
    The resulting high-dimensional densities are discretized on tensor-product grids and represented in tensor-train format, exploiting Kronecker structure for scalable computation. 
    Preconditioning based on estimated posterior moments further improves tractability. 
    We test the method with numerical experiments on EUV reflectometry data for a ruthenium--silicon ($\ce{Ru}$ on $\ce{Si}$) specimen and compare the learned tensor-train posterior to ensemble Markov chain Monte Carlo (MCMC) reference samples. 
    We demonstrate a use case in which  we are able to reuse the tensor-train posterior to obtain samples for shifted target distributions with less computational effort than a comparable MCMC method.
\end{abstract}

\section{Introduction}
    For many applications in high-precision manufacturing in the semiconductor industry, it is essential to have a good understanding of the optical properties of the involved materials. Measuring these properties in the extreme ultraviolet (EUV) range, at wavelengths around \qty{13.5}{nm}, is of high interest for the development of components such as mirror components and photomasks \cite{ciesielski2022forwardmodel}. 
    Since direct measurement of the optical properties of these complex, usually multilayered, systems is generally not possible, reflection measurements \cite{scholze2001highaccuracyeuv,Saadeh2023bayesianinferencestime} can be used instead. 
    Obtaining the parameters of interest from such indirect measurements is mathematically formalized by an \emph{inverse problem}. 
    In previous works \cite{heidenreich2015b, tarantola1987inverse, henze2023bayesianprioranalysis}, this type of problem has been found to be ill-posed, i.e., the solution might not be unique or its dependence on the data might not be continuous, making the parameter reconstruction unreliable in the presence of measurement noise. 
    The inverse problem for reflectometry has been investigated with various methods, e.g., least squares \cite{ciesielski2022forwardmodel} or maximum likelihood estimation \cite{Henn_scatterometry_2012}.
    To improve the quality and reliability of the results, a common technique is to measure the same specimen at a variety of incident angles and photon energy levels, which can also be described by their corresponding wavelengths, thus increasing the data density.

    Instead of seeking a single best reconstruction, one can fix a prior distribution on the parameters and then update it by conditioning it on the measurements.
    The resulting distribution is given by Bayes' theorem.
    In this approach, known as the~\emph{Bayesian inverse problem}, the parameters, as well as quantities of interest that depend on them, can be estimated together with their uncertainty.
    The downside is the high dimensionality and computational cost of the related numerical algorithms.
    A typical method for solution of Bayesian inverse problems is to obtain an approximate posterior sample with Markov chain Monte Carlo (MCMC) simulations \cite{Saadeh_rutheniumsample_2021}.

    While the classical Bayesian MCMC approach is a well-established method for solving this type of inverse problem, it can be computationally expensive, especially when the underlying forward model is expensive to evaluate. Furthermore, it has to be performed anew for each set of measurement data, rendering it inefficient for large numbers of measurements that show only minor shifts in the actual parameters of interest. Such a setting is common, for example, in manufacturing processes, where specimens are produced to known target specifications and ideally exhibit only minimal deviations from them. 
    In case the direct measurement of the quantity of interest is not possible, the necessity of solving an inverse problem can make process supervision challenging.

\begin{figure}[hb]
        \centering
        \resizebox{\linewidth}{!}{\input{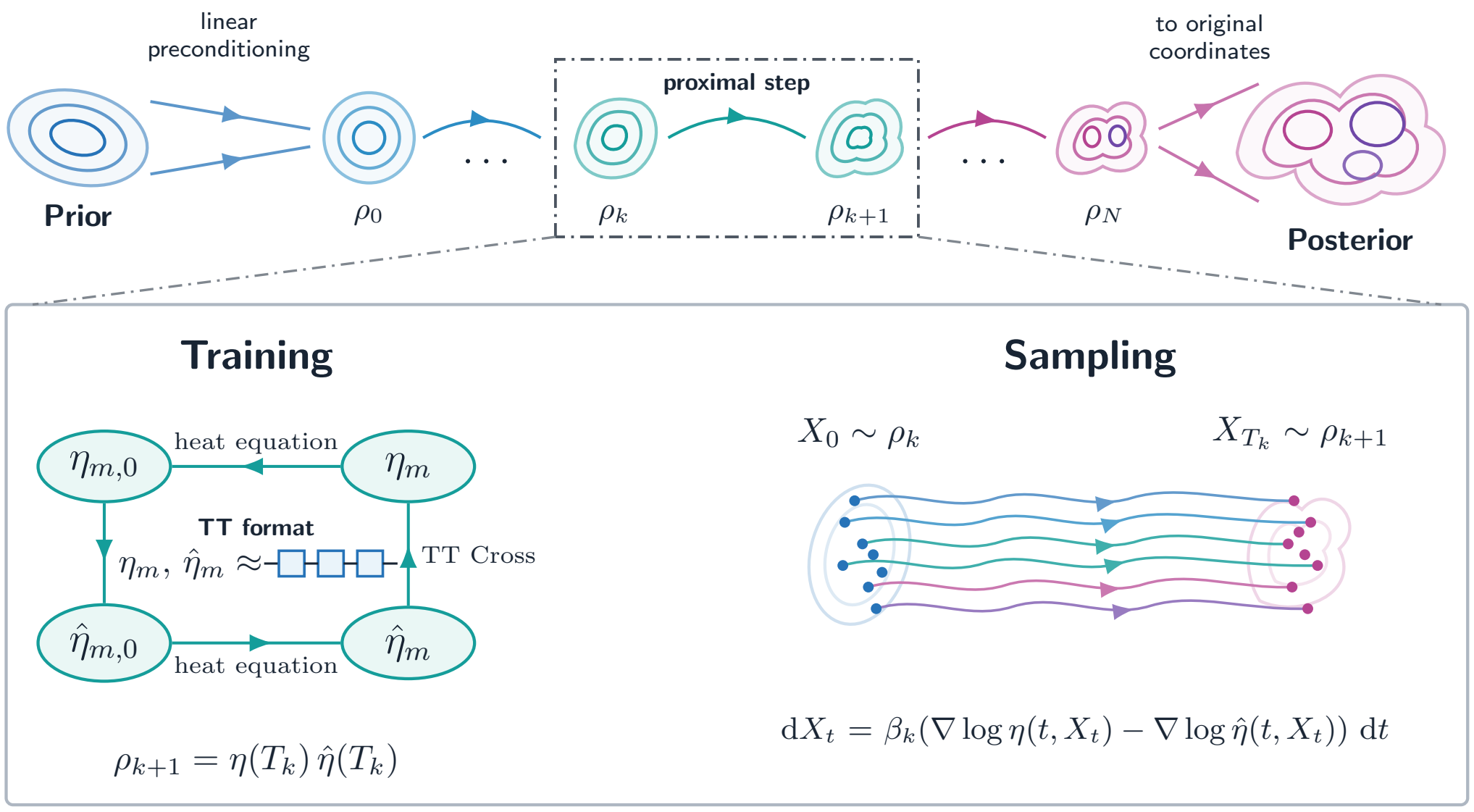}}
        \caption[Conceptual scheme of the rJKOtt method]{Conceptual scheme of the rJKOtt method. Top: the prior is mapped by the linear preconditioner $P$ (\autoref{sec:preconditioning}) to the initial density $\rho^0$; a sequence of regularized JKO (proximal) steps $\rho^k \to \rho^{k+1}$ yields $\rho^{N_{\text{JKO}}}$, which is mapped back to the original coordinates to approximate the posterior. Bottom left (training, \autoref{alg:rjko-proximal}): within one proximal step, the potentials $\eta_m$ and $\hat\eta_m$ are stored in TT format and updated by a fixed-point iteration that alternates heat-equation solves and TT-cross approximations of the coupling conditions; the next density is $\rho^{k+1} = \eta(T_k)\,\hat\eta(T_k)$. Bottom right (sampling, \autoref{sec:sampling}): particles $X_0 \sim \rho^k$ are transported by the ODE~\eqref{eq:ode_interpolation_dynamic} to $X_{T_k} \sim \rho^{k+1}$.}
        \label{fig:rJKOtt_main_illustration}
    \end{figure}

    In the present paper, we consider the regularized Jordan--Kinderlehrer--Otto (JKO) scheme in tensor-train format (rJKOtt), an approach to Bayesian inverse problems that was recently introduced in~\cite{aksenov2025eulerian}.
    The method starts from a tractable distribution  and produces successive approximations of the posterior using a proximal descent in the Wasserstein space, i.e., the space of probability measures endowed with the geometry induced by optimal transport~\cite{ambrosio2005gradient}.
    It is schematically depicted in~\autoref{fig:rJKOtt_main_illustration}.
    At each step, we solve both for the intermediate densities and for the stochastic process that allows one to interpolate between them. 
    The transition between the steps is described by a system of partial differential equations, which can be efficiently solved using the tensor-train (TT) decomposition. 

    The resulting TT approximation can serve as a surrogate model for the posterior and be used to approximately sample from it.
    We demonstrate that, if combined with a suitable preconditioning strategy, the method can produce the posterior approximation using a smaller computational effort than the baseline MCMC.
    Moreover, the approximation by our method can be re-used to provide a warm-start for a solution for a similar measurement, further improving the computational efficiency.

    The rest of the paper is organized as follows.
    In~\autoref{sec:forward}, the forward model is presented in detail.
    The overview of the mathematical concepts used in our method for the inverse problem is given in~\autoref{sec:rJKOtt}, while the numerical scheme as well as the implementation details are described in~\autoref{sec:tensor_train}. 
    Our numerical evaluation of the method for particular reflectometry experiments is given in~\autoref{sec:numerical_experiments}.
    
\section{Forward model}\label{sec:forward}
\begin{figure}[ht]
    \centering
    \includegraphics[width=0.75\linewidth]{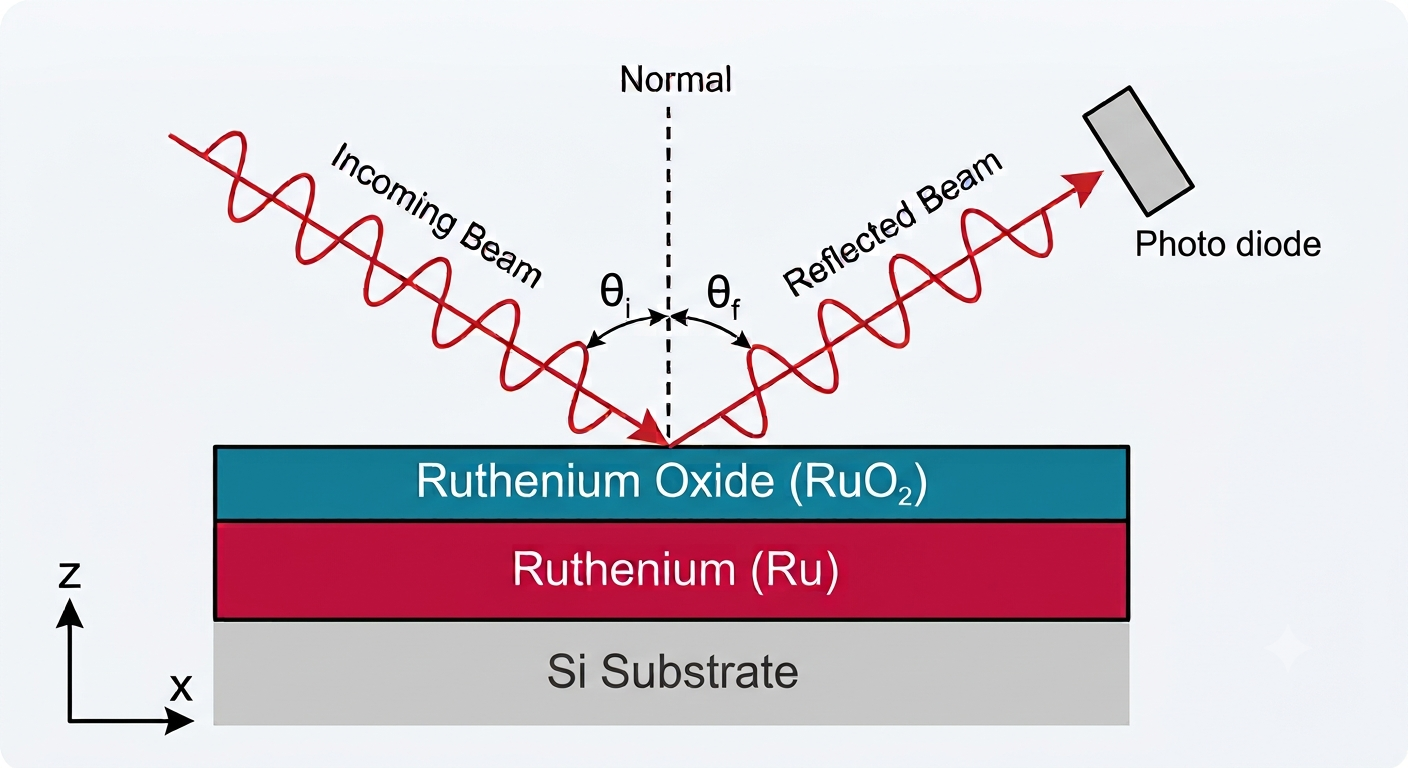}
    \caption{Schematic depiction of the experimental setup and the layer system. The layer thicknesses are not at proper scale w.r.t. each other.}
    \label{fig:layer_system_schematic}    
\end{figure}

    Solving the inverse problem using MCMC requires a suitable forward model to compute the required likelihood function. Since MCMC usually requires a large amount of samples and thus forward model evaluations, this model needs to be reasonably efficient and accurate. Our measurement data contain  the absolute intensity of the measured reflected signal, so our forward model has to provide a reliable mapping between the reflection intensity  and our actual quantities of interest: the optical constants $n, k$, which encode the optical response of our sample, and thickness $h_j$ of the material. The reflection intensity is usually described by ratio of the measured intensity over the incident intensity for any given set of parameters $x$, called the reflectance $R(x)$.
    For the layered systems, we use the Transfer Matrix Method (TMM) described by \cite{ciesielski2022forwardmodel, daillant2008matrixmethod}, to model the transmission and reflection of light through the layered specimen (cf.~\autoref{fig:layer_system_schematic}). 

    Let $j = 0, \dots, N_{L}$ be the layers in the layered system, with $0$ indexing the vacuum layer at the top and $N_{L}$ denoting the substrate layer, which can be considered infinitely thick for computational purposes. The individual layers are considered to be homogeneous and each consist of an individual compound. 
    The light traveling through the system is considered as a monochromatic plane wave $u ^{\pm}(z) \thicksim \exp(\pm i k_{z,j} z)$, with $z$ denoting the height coordinate and $k_{z,j} = (n_j+ik_j)k_{\mathrm{vacuum}}$ being the wave vector component of a wave in layer $j$ in the $z$-direction. All measurements used were performed using only s-polarized light \cite{hansen2025atmoceuv}, whose electric field vector oscillates perpendicular to the incidence plane. Thus, only light with this type of polarization has to be accounted for in the forward model. For each layer, the translation matrix
    \begin{align*}
    \mathcal{T}_{j}(x) = \begin{bmatrix}
      \exp(ik_{z,j} h_j) & 0\\
      0 & \exp(-ik_{z,j}h_j)\\
      \end{bmatrix}
      \end{align*}
      describes the behavior of the wave traveling through the material while the refraction matrix

\begin{align*}
     \mathcal{R}_{j-1,j}(x) = \left[{\begin{array}{*{20}{c}}{{p_{j/j - 1}}}&{{m_{j/j - 1}}}\\{{m_{j/j - 1}}}&{{p_{j/j - 1}}}\end{array}} \right],
\end{align*}
with 
      \begin{align*}
          \begin{split}{p_{j/j - 1}} &= \frac{{{k_{z,j}} + {k_{z,j - 1}}}}{{2{k_{z,j}}}} \cdot \exp \left({- \frac{1}{2}{{({k_{z,j}} + {k_{z,j - 1}})}^2} \cdot s_{j/j - 1}^2} \right), \\ {m_{j/j - 1}} &= \frac{{{k_{z,j}} - {k_{z,j - 1}}}}{{2{k_{z,j}}}} \cdot \exp \left({- \frac{1}{2}{{({k_{z,j}} - {k_{z,j - 1}})}^2} \cdot s_{j/j - 1}^2} \right).\end{split}
      \end{align*}
      
      describes their behavior at the layer interfaces for a given set of optical constants (OC) and geometric parameters, all contained in $x$.  
      Additionally, a roughness term $\exp \left({- \frac{1}{2}{{({k_{z,j}} \pm {k_{z,j - 1}})}^2} \cdot s_{j/j - 1}^2} \right)$ can  be used to model the effects of rough surfaces and other imperfections at the layer interfaces, which usually results in an additional attenuation of the signal at the layer interface transitions. The derivation is described in detail in \cite{daillant2008matrixmethod}.
    This yields a forward model that is fast to evaluate, as well easily extendable to multiple layers, with each layer just adding their respective refraction matrices $\mathcal{R}_{j-1,j}(x)$ and translation matrices $\mathcal{T}_{j}(x)$ to the system of equations 

    \begin{align}
      \begin{bmatrix}
      u^+(z_0) \\
      u^-(z_0) \\
      \end{bmatrix} = M(x) \cdot 
      \begin{bmatrix}
      u^+(z_{N_L}) \\
      u^-(z_{N_L}) \\    
      \end{bmatrix} =  \begin{bmatrix}
            M_{11}(x) & M_{12}(x) \\
            M_{21}(x) & M_{22}(x)
        \end{bmatrix}
        \cdot 
      \begin{bmatrix}
      u^+(z_{N_L}) \\
      u^-(z_{N_L}) \\    
      \end{bmatrix} , 
      \label{eq:matrixmethod}
    \end{align}   
    where 
    \begin{align}
        M(x) = \mathcal{R}_{0,1}(x) \mathcal{T}_{1}(x) \mathcal{R}_{1,2}(x) \mathcal{T}_{2}(x) \cdots  \mathcal{T}_{N_L-1}(x) \mathcal{R}_{N_L-1,N_L}(x) .
    \end{align}

   Since the thickness of the substrate is considered to be large compared to the wavelengths of the used light source as well as compared to the thin material film, no upward returning wave $u^+(z_n)$ from the bulk material is to be expected. The reflectance coefficient $R(x)$ can then be computed as the ratio of the intensities of the reflected and incident waves at the vacuum interface,
    \begin{align}
        R(x) = \left |\frac{u(z_0)^+}{u(z_0)^-}\right |^2 = \left |\frac{M_{12}(x)}{M_{22}(x)} \right | ^2 \text{.}
    \end{align}

\section{Regularized JKO scheme} \label{sec:rJKOtt}
    In this section, we describe the mathematical principles behind our method.
    We rely on the theory of Wasserstein spaces.
    The Wasserstein distance is based on the solution of the optimal transport problem:
    \begin{equation}
        W^2_2(\rho_0, \rho_1) := \inf_{\gamma \in \Gamma(\rho_0,\rho_1)} \int \|x - y\|^2_2 \mathrm{d}\gamma(x, y),
    \end{equation}
    where $\Gamma(\rho_0, \rho_1)$ is the set of all \emph{transport plans}, i.e., joint probability distributions of $(x, y)$ with first and second marginals $\rho_0$ and $\rho_1$ respectively.
    The rich geometric structure of this space~\cite{ambrosio2005gradient,villani2021topics} allows one to adapt numerical minimization methods to various problems, including the posterior approximation in Bayesian inverse problems.
    In the following, we explain the formulation of our method in this space.
    \subsection{Minimization on the space of probability measures}
        Our method fits into the broader approach of Variational Inference~\cite{bishop2006pattern,yao2022mean}.
        In this approach, a tractable model of the target distribution $\rho_\infty$ is constructed by minimizing a functional that quantifies the discrepancy between $\rho_\infty$ and the current approximation.
        A prominent example of such functionals is the Kullback--Leibler divergence~\cite{kullback1951information,mohammad2015entropy}:
        \begin{equation}\label{eq:KL}
            \operatorname{KL}(\rho | \rho_\infty) = 
                \begin{cases}
                    \int\log\left( \frac{\mathrm{d}\rho}{\mathrm{d}\rho_\infty}\right)\mathrm{d}\rho, \quad \text{if } \rho \ll \rho_\infty \\
                    +\infty,\quad \text{otherwise.}
                \end{cases}
        \end{equation}
        Its minimizer does not change if the posterior distribution is replaced by its unnormalized estimate: $\log \rho_\infty \to \log\rho_\infty + \text{const}$, making it particularly relevant in the Bayesian inversion case~\cite{blei2017variational}, where the normalizing constant is not directly accessible.
    
        To minimize the KL divergence, we rely on the Jordan--Kinderlehrer--Otto (JKO) scheme, introduced in the seminal paper~\cite{jordan1998variational}
        \begin{equation}\label{eq:JKO}
            \rho^{k+1} \in \underset{\rho}{\arg\min} \operatorname{KL}(\rho|\rho_\infty) + \frac{1}{2T_k}W^2_2(\rho, \rho^k).
        \end{equation}
        Here $\rho^k$ is the sequence, generated by the method, $W_2$ denotes the Wasserstein distance and $T_k$ is a step size parameter.
        Essentially, this is a proximal method~\cite{combettes2011proximal}, where at each step, one attempts to minimize the target functional at the next iterate, while not diverging too far from the previous iterate.
        Small values of the step size parameter $T_k$ cause larger distance penalty and vice versa.

        In the construction of our method, we follow the ``optimize-then-discretize'' paradigm~\cite{hinze2008optimization}. 
        The next iterate of the minimization method is defined as a solution of a system of partial differential equations.
        A discrete approximation of this solution is then computed, using a suitable discretization method.
        
    \subsection{Representation as a coupled PDE system}
        To be able to numerically approximate the JKO step using TT discretization, a few adjustments have to be made.
        First, we switch to the dynamic formulation of the Wasserstein distance, due to Benamou and Brenier~\cite{benamou2000computational}.
        Second, entropic regularization is added to the dynamics.
        We arrive at a PDE-constrained optimization problem:
        \begin{align}\label{eq:reg_dyn_JKO}
            \rho^{k+1} =  &\arg\min\limits_{v,\rho }\int\limits_0^{T_k}\int\limits_{\mathbb{R}^d} {\frac{1}{2}\|v(t,x)\|^2}{\rho(t, x)}\diff{x}\diff{t} + \operatorname{KL}(\rho({T_k},\cdot) | \rho_\infty), \\
                    &\operatorname{s.t.}\quad   
                    \begin{aligned}\notag
                        \partial_t \rho + \nabla\cdot(\rho v) &= \beta_k\operatorname{\Delta}\rho,\\
                        \rho(0, x) &= \rho^{k}(x).
                    \end{aligned}
        \end{align}
        In the Benamou--Brenier formulation, a vector field $v$ is introduced, that smoothly interpolates the distribution in time $t \in [0, T_k]$, which is enforced by the continuity equation $\partial_t \rho + \nabla\cdot(\rho v) = 0$.
        The kinetic energy term in~\eqref{eq:reg_dyn_JKO} quantifies the transportation cost along $v$.
        The entropic regularization $\beta_k \Delta\rho$ is commonly added into the continuity equation~\cite{cuturi2013sinkhorn,peyre2015entropic,carlier2017convergence} in order to make the problem strongly convex.

        Considering the first-order criticality conditions of~\eqref{eq:reg_dyn_JKO} results in a system consisting of a Fokker--Planck equation and a Hamilton--Jacobi--Bellman equation with nonlinear coupling both in the equations, initial and terminal conditions.
        Using the Hopf--Cole transform, we introduce new variables 
        \begin{equation}\label{eq:Hopf-Cole}
            \eta(t,x) := e^{\frac{\Phi(t,x)}{2\beta_k}} \text{ and }
            \hat{\eta}(t,x) := \rho(t,x)e^{\frac{-\Phi(t,x)}{2\beta_k}} 
        \end{equation}
        where the optimal vector field  has a gradient structure: $v(t, x) = \nabla \Phi(t, x)$.
        In the new variables, the system takes the form
        \begin{gather}\label{eq:heat_system}
            \partial_t \hat{\eta} = \beta_k\Delta\hat{\eta} \\
            \partial_t {\eta} = -\beta_k\Delta\eta \label{eq:heat_direct}\\
            \hat{\eta}(0, x) = \frac{\rho^k(x)}{\eta(0, x)} \label{eq:fp_initial_KL}\\ 
            \eta(T_k, x) = \left(\frac{\rho_\infty}{\hat\eta(T_k, x)} \right)^{\frac{1}{1 + 2\beta_k}} .\label{eq:fp_terminal_KL}
        \end{gather}
        The detailed derivation can be found in~\cite{li2020fisher,li2023kernel}.
        We note that the terminal condition is specific for the target functional $\operatorname{KL}(\cdot | \rho_\infty)$.
        % $\eta(T_k,x) = e^{-\frac1{2\beta_k}{\delta\mathcal{E}(\eta(T_k, \cdot), \hat\eta(T_k, \cdot), x)}}$

        The system {\eqref{eq:heat_system}-\eqref{eq:fp_terminal_KL}} is more favorable to solve than the original optimality system for several reasons.
        First, the coupling is now present only in the initial and terminal conditions.
        Thus, only one initial-boundary value problem has to be solved at a time.
        Satisfying the coupled initial and terminal conditions can be reformulated as a fixed-point problem and efficiently solved with accelerated numerical methods, see~\cite{bungartz2016precice,aksenov2021application,delaisse2023quasi}.
        Second, as the Laplace operator has a Kronecker sum structure~\cite{graham2018kronecker}, the solution of the heat equation has a Kronecker product structure.
        For example, for~\eqref{eq:heat_system},
        \begin{equation}
            \hat\eta(t, \cdot) = \exp\left({\beta_k t \frac{\partial^2}{\partial x_1^2}}\right) \otimes \dots \otimes \exp\left({\beta_k t \frac{\partial^2}{\partial x_d^2}}\right) \hat\eta(0, \cdot).
        \end{equation}
        Essentially, the heat equation in different dimensions is decoupled.
        In the sequel, we show how this structure allows one to efficiently solve the heat equation in the TT format.
        
        \subsection{Sampling dynamics}\label{sec:sampling}
        The variables $\eta,\ \hat\eta$ encode the information on both the dynamically changing density and the vector field:
        \begin{equation}
            \rho(t, x) = \eta(t, x)\hat{\eta}(t, x),\quad  v(t, x) = 2 \beta_k \nabla_x \log{\eta(t, x)}.
        \end{equation}    
        Solving the aforementioned equations gives us the \emph{Eulerian perspective} on the dynamics, i.e., the transient change of functions in the fixed frame of reference $x$.
        A passage can be made to \emph{Lagrangian perspective}, where the trajectories of individual moving particles are considered.
        In particular, the Fokker--Planck equation
        \begin{equation}
            \partial_t \rho + 2\beta_k\nabla\cdot(\rho \nabla \log\eta)  = \beta_k \Delta \rho
        \end{equation}
        corresponds to the dynamics governed by the following SDE~\cite{pavliotis2014stochastic}:
        \begin{equation}\label{eq:sde_interpolation_dynamic}
            \mathrm{d}X_t = 2\beta_k\nabla{\log\eta(t, X_t)}\mathrm{d}t + \sqrt{2\beta_k} \mathrm{d}W_t,\qquad
            X_0 \sim \rho^k.
        \end{equation}
        After transforming the Fokker--Planck equation  as follows
        \begin{multline*}
        \partial_t \rho + 2\beta_k\nabla\cdot(\rho \nabla \log\eta) - \beta_k \Delta \rho \\ 
            =  \partial_t \rho + 2\beta_k\nabla\cdot(\rho \nabla \log\eta) - \beta_k \nabla \cdot (\nabla \rho)
            =  \partial_t \rho + 2\beta_k\nabla\cdot(\rho \nabla \log\eta) - \beta_k \nabla \cdot (\rho \nabla \log \rho) \\
            =  \partial_t \rho + \beta_k\nabla\cdot\left(\rho \nabla (2\log\eta - \log \rho)\right) 
            =  \partial_t \rho + \beta_k\nabla\cdot\left(\rho \nabla (\log\eta - \log \hat\eta)\right) 
        \end{multline*} 
        it becomes a continuity equation with resulting vector field $\tilde v(t, x) = \beta_k \nabla_x \left(\log \eta(t, x) - \log\hat\eta(t, x)\right)$.
        This corresponds to the dynamics governed by the ODE~\cite[{Lemma {8.1.4}}]{ambrosio2005gradient}
        \begin{equation}\label{eq:ode_interpolation_dynamic}
            \diff X_t = \beta_k \nabla \left(\log \eta(t, X_t) - \log\hat\eta(t, X_t)\right)\diff t.
        \end{equation}
        If the current distribution $\rho^k$ can be efficiently sampled from, and the solution $(\eta,\ \hat\eta)$ is known, then solving either the SDE~\eqref{eq:sde_interpolation_dynamic} or the ODE~\eqref{eq:ode_interpolation_dynamic} with initial values distributed as $\rho^k$ gives a sample, approximately distributed as $\rho^{k+1}$.
        \subsection{Convergence}
        \begin{corollary}
            Under the regularity assumptions on the moments and the entropies of the generated sequence, stated in  \cite[{Corollary 4}]{aksenov2025eulerian}, if 
            \[
                -\nabla^2 \log\rho_\infty(x) \succeq \lambda\Id\quad\forall x \in \mathbb{R}^d
            \] 
            for some $\lambda > 0$, after $k$ steps of the regularized JKO scheme with constant step size $T_k \equiv T > 1/\lambda$ and regularization $\beta$, the error can be quantified as
            \[
                W^2_2(\rho^{k}, \rho_\infty) \leq \left( \frac{2}{1 + T\lambda}\right)^k W^2_2(\rho^0, \rho_\infty)
                + \frac{8K^*}{\lambda - \sfrac{1}{T}}\beta \log\left(\frac{1}{\beta} \right)
            \] for some $K^* > 0$.
        \end{corollary}
        The first term describes the geometric contraction of the initial error with the number of steps, while the second term is the bias introduced by the entropic regularization, which vanishes as $\beta \to 0$. Since the JKO scheme is an implicit method, arbitrarily large step sizes $T_k$ can be used.
        The tradeoff is that the fixed-point iteration converges more slowly for larger $T_k$.
        We show numerically in~\autoref{sec:numerical_experiments} that, if fixed-point acceleration techniques~\cite{tang2022accelerating,walker2011anderson} are used, the problem can still be efficiently solved in one or two proximal steps. 
        
\section{Numerical implementation}\label{sec:tensor_train}
    The main difficulty in the Eulerian approach is that the number of degrees of freedom in the discretization scales exponentially with the problem dimension.
    In case of optical inverse problems, the dimension is the number of the specimen's parameters to reconstruct.
    To mitigate this, we compress the discretization coefficients using the TT decomposition.
    It is also important that the discretization grid is built in the region of high posterior concentration.
    Thus, we introduce preconditioning in form of coordinate transformation.
    The outline of the complete method is given in~\autoref{alg:rjko-main} at the end of this section.
    \subsection{Tensor-train discretization}
        The problem's dimension $d$ is defined by the number of parameters, the joint posterior distribution of which we want to infer.
        A naive discretization of~{\eqref{eq:heat_system}-\eqref{eq:fp_terminal_KL}} with e.g., Finite Element or Finite Difference methods would result in a system with $\mathcal{O}(N^d)$ degrees of freedom.
        Even for a moderate dimension of $d \geq 4$, storing and processing such amounts of data becomes infeasible in practice. 
        This is known as the \emph{curse of dimensionality}.
        To mitigate that, we employ low-rank compression techniques, in particular, the TT decomposition.

        The TT decomposition, introduced in~\cite{oseledets2011tensor,hubener2010concatenated}, can be viewed as a multidimensional generalization of the truncated SVD~\cite{oseledets2009breaking}.
        A multidimensional tensor $g \in \mathbb{R}^{N_1 \times \dots \times N_d}$ is said to be in the TT format if it can be represented as
        \begin{equation}\label{eq:tt_definition}
            g_{i_1\dots i_d} = G^1_{i_1 l_1}G^2_{l_1 i_2 l_2}\dots G^d_{l_{d-1}i_d},
        \end{equation}
        where for each $n = 1, \dots, d$, $G^n_{l_{n-1}i_n l_n} \in \mathbb{R}^{r_{n - 1} \times N_n \times r_n}$ are $3$-dimensional tensors, called \emph{tensor cores}. 
        We assume $r_0 = r_{d} = 1$ and Einstein summation convention.
        The values $r_k$ are called \emph{TT-ranks}.
        Let us denote by $r := \max_k r_k$ the TT rank of $g$, and let $N := \max_k N_k$.

        The storage of a tensor in the TT format requires $\mathcal{O}(dNr^2)$ memory, making it tractable for larger dimensions, provided that the rank stays low.
        The access to a given element requires $\mathcal{O}(dr^2)$ operations, unlike $\mathcal{O}(1)$ in case of storing a full tensor.
        Fortunately, linear algebra operations, such as summation or contraction of indices, can also be efficiently implemented in the TT format, i.e., without computing each of the tensors' elements, with complexity polynomial in $d$.
        Moreover, there exist tensor reconstruction algorithms, such as cross-approximation~\cite{oseledets2010tt}, that allow one to acquire a tensor-train approximation of a dense tensor, by computing only a limited amount of its values.
        Finally, the TT rounding algorithm~\cite{oseledets2011tensor} can be used to reduce the TT rank, while controlling the approximation error.
        To conclude, the diverse toolbox of TT algorithms allows one to solve various problems in the compressed format. For additional reference, see~\cite{phan2016tensor,cichocki2016tensor}.

        In our implementation of the regularized JKO scheme, at $k$-th iteration, the functions $\eta^k, \hat\eta^k$ are discretized in space in a finite-difference fashion, i.e., we look for tensors, representing their values at a predefined rectangular grid:
        \[
            \eta^k_{i_1, \dots, i_d}(t) := \eta^k(t, x_{1, i_1}, \dots, x_{d, i_d}),
        \]
        where the grid in each dimension is in interval $[L_n,\ R_n]$, and the grid nodes are defined as
        \[
            x_{n, i} = L_n + h_n\cdot i ,\qquad h_n = \frac{R_n - L_n}{N_n}.
        \]
        We note that the spatial discretization in a finite-element fashion, i.e., working with a tensor of coefficients of functions' representation in some basis, is also possible in the framework of functional tensor trains~\cite{oseledets2013constructive,gorodetsky2019continuous}.

        Using the Kronecker product structure of the solution of the heat equation, we can efficiently solve it in the TT format while preserving the TT rank.
        Let us denote by $L_h$ the finite dimensional approximation of the Laplace operator:
        \begin{equation}
            \label{eq:discrete_Lap}
            L_h = D_{1, h} \otimes I_2 \otimes \cdots \otimes I_d + 
                I_1 \otimes D_{2, h} \otimes \cdots \otimes I_d + \dots + 
                I_1 \otimes I_2 \otimes \cdots \otimes D_{d, h},
        \end{equation}
        where $\otimes$ denotes the Kronecker product, $I_k \in \mathbb{R}^{N_k \times N_k}$ is an identity matrix and 
        \begin{equation*}
            D_{n, h} = \frac{1}{h_n^2}
                \begin{pmatrix}
                    -1 & 1      &           &        &         \\
                    1  & -2     & 1         &        &         \\
                       & \ddots & \ddots    & \ddots &         \\
                       &        & \ddots    & \ddots & 1       \\
                       &        &           & 1      & -1        
                \end{pmatrix}.
        \end{equation*}
        The solution operator $\mathcal{H}_d(\beta_k t)$ of the heat equation at time $t$ has the structure
        \begin{equation}
            \hat\eta(t) = \mathcal{H}_d(\beta_k t)\hat\eta(0) := e^{\beta_k t L_h}\hat\eta(0).
        \end{equation}
        As noted in, for example~\cite{chertkov2021solution}, we can write for~\eqref{eq:heat_system},
        \begin{multline*}
            \hat\eta(t)_{i_1\dots i_d} 
            = \left(e^{\beta_k t L_h} \right)_{i_1\dots i_d j_1 \dots j_d}\hat\eta(0)_{j_1 \dots j_d} \\
            = \left(e^{\beta_k t L_h} \right)_{i_1\dots i_d j_1 \dots j_d} \hat\eta^1_{j_1 l_1}\hat\eta^2_{l_1 j_2 l_2}\dots \hat\eta^d_{l_{d-1}j_d} 
            = \left(e^{\beta_k t D_{1, h}}_{i_1 j_1} \dots e^{\beta_k t D_{d, h}}_{i_d j_d}\right) \hat\eta^1_{j_1 l_1}\hat\eta^2_{l_1 j_2 l_2}\dots \hat\eta^d_{l_{d-1}j_d} \\
            = \left(e^{\beta_k t D_{1, h}}_{i_1 j_1} \hat\eta^1_{j_1 l_1}\right) \left(e^{\beta_k t D_{2, h}}_{i_2 j_2} \hat\eta^2_{l_1 j_2 l_2} \right)\dots \left(e^{\beta_k t D_{d, h}}_{i_d j_d} \hat\eta^d_{l_{d-1}j_d}\right).
        \end{multline*}
        where $\hat\eta^n$ are the cores of TT representation of $\hat\eta(0)$, and $e^{\beta_k t D_{n, h}} \in \mathbb{R}^{N_n \times N_n}$ denotes the matrix exponential.
        One can see that each matrix exponential is contracted with the respective tensor core.
        Thus, the solutions of the heat equation in the individual dimensions are decoupled. 
        For fixed grid and value of $\beta_k t$, the matrix exponentials can be precomputed, and thus the solution of the heat equation can be performed in $\mathcal{O}(dN^2r^2)$.
        Note that $T_k$ enters the computational formulas only together with $\beta_k$.
        Thus it is more convenient in the sequel to refer to $T_k\beta_k$ as effective time parameter.

        Since the time evolution of $\hat\eta$ from $0$ to $T_k$, and of $\eta$ in the inverse direction, can be performed effectively without time marching, we only need to find approximations for $\eta$ and $\hat\eta$ at initial or terminal time, satisfying the coupling conditions.
        If one of the boundary parameters is available in the TT format, the other one can be reconstructed with cross-approximation using~\eqref{eq:fp_initial_KL} or~\eqref{eq:fp_terminal_KL}.
        To satisfy both conditions at once, we employ the fixed-point iteration, illustrated in the bottom-left panel of~\autoref{fig:rJKOtt_main_illustration}.
        The pseudocode for one step is given in~\autoref{alg:rjko-proximal}.
        \begin{figure}[hb]
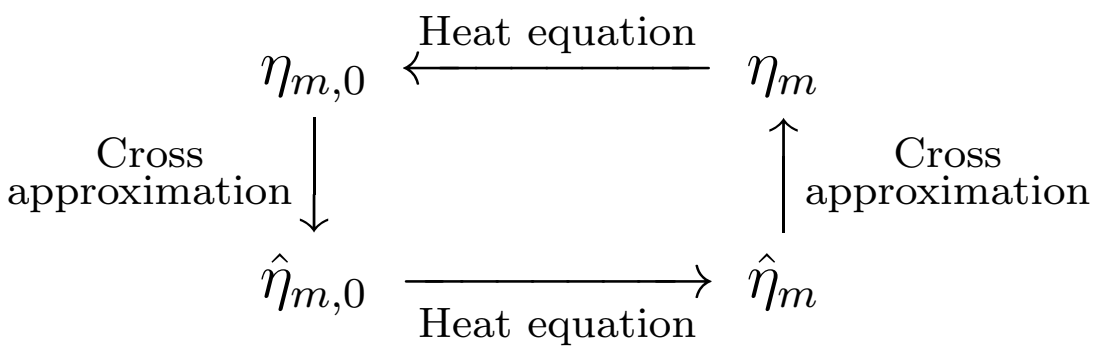

            \centering
            \[ \large
            \begin{CD}
                \eta_{m, 0}           @<{\text{Heat equation}}<<    \eta_{m} \\
                @V{\substack{\text{Cross} \\ \text{approximation}}}VV             @AA{\substack{\text{Cross} \\ \text{approximation}}}A \\   
                \hat\eta_{m, 0}        @>>{\text{Heat equation}}>     \hat\eta_{m}
            \end{CD}
            \]
            \caption{Scheme of the fixed-point iteration}
            \label{fig:fp_cycle}
        \end{figure}  
    \begin{algorithm}[ht]
        \caption{One regularized JKO proximal step in TT format}
        \label{alg:rjko-proximal}
        \small      
        \begin{algorithmic}[1]
            \Require TT initial guess $\eta_0$; fixed-point tolerance
            $\varepsilon_{\mathrm{FP}}$; positivity floor $\varepsilon_{\mathrm{pos}}$;
            TT~tolerance $\varepsilon_{\mathrm{TT}}$; maximum TT~rank $r_{\max}$;
            maximum iterations $M_{\max}$; fixed-point update rule
            \Ensure converged potentials and reconstructed density $\rho^{k+1}$
            \Function{TTProximalStep}{$\rho^k,\rho_\infty,T_k,\beta_k$}
            
            \State $m\gets0$
            \Repeat     
              \State $\eta_{m,0}\gets\Heat(\beta_kT_k)\eta_m$
                \Comment{Backward heat solve}
              \State $\displaystyle
              \widehat\eta_{m,0}\gets\TTcross\!\left[
              z\mapsto\frac{\rho^k(z)}
              {\max\{\eta_{m,0}(z),\varepsilon_{\mathrm{pos}}\}}
                \right]$      
              \State $\widehat\eta_m\gets\Heat(\beta_kT_k)\widehat\eta_{m,0}$
                \Comment{Forward heat solve}
              \State $\displaystyle
              \widetilde\eta_m\gets\TTcross\!\left[
              z\mapsto\left(
              \frac{\rho_\infty(z)}
              {\max\{\widehat\eta_m(z),\varepsilon_{\mathrm{pos}}\}}
              \right)^{1/(1+2\beta_k)}
            \right]$      
            \State        
              % $g_m\gets\widetilde\eta_m$;
              $r_m\gets \widetilde\eta_m-\eta_m$;
              $e_m\gets\RelRes{r_m}{\eta_m}$
            
              \State $\displaystyle
              \eta_{m+1}\gets
              \Call{FixedPointStep}{\{\eta_j,\widetilde\eta_j,r_j\}_{j=0}^{m}}$
              \Comment{For example, Picard or Anderson}
              \State $\eta_{m+1}\gets
                \TTround(\eta_{m+1}, \varepsilon_{\mathrm{TT}},r_{\max})$
              \Comment{Preserve low rank}
              \State $m\gets m+1$
            \Until{converged or $m = M_{\max}$}
            
            \State $\rho^{k+1}\gets
              \operatorname{Normalize}(\eta_M\widehat\eta_M)$
            \State \Return
            $\eta_M,\widehat\eta_{M,0},\widehat\eta_M,\rho^{k+1}$
            \EndFunction        
        \end{algorithmic}
    \end{algorithm}

    Unlike particle-based samplers such as MCMC, Stein variational gradient descent (SVGD) or Langevin dynamics, our method evaluates the posterior only at predefined grid points, which allows one to implement a caching strategy.
    In practice, this allows one to significantly reduce the number of costly posterior evaluations.
    The TT solution can be used to generate samples from the approximate posterior by numerically simulating the dynamics~\eqref{eq:sde_interpolation_dynamic} or \eqref{eq:ode_interpolation_dynamic}.
    Alternatively, the next density $\rho^{k+1} = \eta(T_k) \hat\eta(T_k)$ can be assembled as a TT of rank at most $r^2$.
    It can then be used to sample with the conditional marginal method~\cite{dolgov2020approximation}.
    An arbitrary number of independent samples can be generated without additional posterior calls, resulting in a lower variance of estimates for quantities of interest.
    This is similar to other training-based methods, such as score~\cite{holzschuh2023solving} or flow matching~\cite{chemseddine2024conditional,kruger2026well}.
    In contrast, for particle-based methods, generating new samples requires posterior calls, and the samples are correlated, both in timestep (autocorrelation) and, for certain methods, within the ensemble. 
    While our method introduces a bias due to regularization, spatial and temporal discretization, we argue that the optimal tradeoff can be achieved for more reliable inference of the parameters.

    \subsection{Preconditioning}\label{sec:preconditioning}
    In order to avoid an unnecessarily large number of grid nodes, it is beneficial that the grid and the initial distribution are centered around the high-probability region of the target distribution. 
    Unfortunately, the prior distribution doesn't always capture this region well due to the concentration of the posterior.
    Additionally, the TT algorithms are observed to yield better low-rank approximation when the target distribution's covariance matrix is close to diagonal. 
    Thus, we precondition the problem by switching to another system of coordinates.
    In particular, we estimate the mean and covariance $(\mu, \Sigma)$ of the posterior, and then use the mapping
    \begin{equation}\label{eq:orig_to_pc_transform}
        z = P(x) := \Sigma^{-\sfrac{1}{2}}(x - \mu)
    \end{equation}
    This transformation is an optimal transport mapping between a Gaussian distribution $\mathcal{N}(\mu, \Sigma)$ and the standard normal distribution.
    We expect that the modified target $P_\sharp \rho_\infty$ is close to the standard normal.
    We then can use the standard normal, which can be efficiently sampled and can be represented as a {rank-$1$} TT, as an initial distribution.
    After running the regularized JKO method and producing samples $\{z_i\}_{i=1}^{N_{\text{samples}}}$, they can be mapped back into the original coordinates:
    \[
        x_i = P^{-1}(z_i) = \Sigma^{\sfrac{1}{2}}z_i + \mu.
    \]
    
    The estimation of posterior mean and covariance can be performed with a combination of a variety of methods, such as Laplace approximation~\cite{scheichl2021numerical} or importance sampling~\cite{agapiou2017importance}.

        \subsubsection{Laplace approximation}
            Laplace approximation is based on matching the second-order Taylor expansion of $\log\rho_\infty$ at its mode with the log-density of a Gaussian distribution:
            \begin{align}
                \mu &= \arg\max_x \log\rho_\infty(x), \label{eq:laplace_MAP} \\
                \Sigma &= \left[ \left. -\nabla^2 \log\rho_\infty(x) \right|_{x=\mu} \right]^{-1}. \label{eq:laplace_covariance}
            \end{align}
            The maximum a posteriori (MAP) estimate of the parameters $\mu$ can be found using derivative-free optimization methods, such as Nelder--Mead~\cite{nelder1965simplex}, or with gradient-based methods, such as \mbox{(L-)BFGS}~\cite{liu1989limited}, trust-region~\cite{conn2000trust} or SQP~\cite{boggs1995sequential}.
            In the latter case, gradients are obtained by automatic differentiation or finite differences.
            As for computing the covariance matrix, if the posterior supports automatic Hessian-Vector products~\cite{pearlmutter1994fast}, they can be used to estimate the inverse Hessian in~\eqref{eq:laplace_covariance} with an iterative matrix inversion method.
            Otherwise, one has to rely on finite-difference approximations of the second derivatives, whose complexity scales quadratically with dimension and thus can become quite costly.
            Both covariance estimation methods can be prone to numerical error, thus, a careful regularization has to be implemented.
        \subsubsection{Importance sampling}
            Importance sampling is a way to estimate the expectation of a function under the posterior by sampling from a tractable distribution $\rho_I$ instead and re-weighting each summand by the density ratio $\rho_\infty/\rho_I$:
            \begin{equation}\label{eq:importance_sampling}
                \mathbb{E}_{x \sim \rho_\infty} f(x) = 
                    \int f(x)\rho_\infty(x) \diff x =  
                    \int f(x)\frac{\rho_\infty(x) }{\rho_I(x)}\rho_I(x) \diff x =
                    \mathbb{E}_{x \sim \rho_I}\left[ \frac{\rho_\infty(x) }{\rho_I(x)}f(x) \right]
            \end{equation}
            where $\rho_I$ is the selected importance distribution.
            For preconditioning, we use this method iteratively, as suggested in~\cite{morzfeld2018iterative}.
            In particular, we start with some initial  $(\mu_0, \Sigma_0)$, obtained, for example, from the prior or from the Laplace approximation.
            Then, we sample from the current Gaussian distribution, evaluate the log density at the new sample, and use the importance sampling formula~\eqref{eq:importance_sampling} to update the mean and covariance value.

    \begin{algorithm}[ht]
        \caption{Preconditioned Eulerian rJKOtt training and sampling}
        \label{alg:rjko-main}
        \small
        \begin{algorithmic}[1]
            \Require tractable prior density $\rho_{\text{prior}}$; unnormalized target
            density $\rho_{\infty}$; affine preconditioner $P(x)=\Sigma^{-\sfrac{1}{2}}(x-\mu)$; outer
            parameters $\{T_k,\beta_k\}_{k=0}^{N_{\text{JKO}}-1}$; TT/grid parameters; number of
            particles $N_{\text{samples}}$
            \Ensure approximate density $\rho^{N_{\text{JKO}}}$ and particles
            $\{X_i^{N_{\text{JKO}}}\}_{i=1}^{N_{\text{samples}}}$ in the original coordinates
            
            \State Define $z=P(x)$ and $P^{-1}(z)=\Sigma^{\sfrac{1}{2}}z+\mu$
            \State $\bar\rho^0=\mathcal N(0,\Id)$
            \State $\bar\rho_\infty(z)\propto
              \rho_\infty(P^{-1}z)$
            \State Draw $Z_i^0\sim\bar\rho^0$ independently for $i=1,\ldots,N_{\text{samples}}$
            
            \For{$k=0,\ldots,N_{\text{JKO}}-1$}
              \State $(\eta_M^k,\widehat\eta_{M,0}^k,\widehat\eta_M^k,
              \bar\rho^{k+1})\gets
              \Call{TTProximalStep}{\bar\rho^k,\bar\rho_\infty,T_k,\beta_k}$
              \State Define, for $t\in[0,T_k]$,
              \Statex \hspace{2.2em}$\displaystyle
              \eta^k(t,\cdot)=\Heat\!\bigl(\beta_k(T_k-t)\bigr)\eta_M^k,$
              \qquad
              $\displaystyle
              \widehat\eta^k(t,\cdot)=\Heat(\beta_k t)\widehat\eta_{M,0}^k$
              \For{$i=1,\ldots,N_{\text{samples}}$ \textbf{in parallel}}
                \State Solve on $[0,T_k]$ with $Z_i(0)=Z_i^k$:
                \Statex \hspace{2.2em}$\displaystyle
                \dot Z_i(t)=\beta_k\nabla_z\!\left[
                \log\eta^k(t,Z_i(t))-\log\widehat\eta^k(t,Z_i(t))
                \right]$
                \State $Z_i^{k+1}\gets Z_i(T_k)$
              \EndFor
            \EndFor
            
            \State $X_i^{N_{\text{JKO}}}\gets P^{-1}(Z_i^{N_{\text{JKO}}})$ for $i=1,\ldots,N_{\text{samples}}$
            \State $\rho^{N_{\text{JKO}}}\gets(P^{-1})_\sharp\bar\rho^{N_{\text{JKO}}}$
            \State \Return $\rho^{N_{\text{JKO}}},\{X_i^{N_{\text{JKO}}}\}_{i=1}^{N_{\text{samples}}}$
        \end{algorithmic}
    \end{algorithm}

\section{Numerical experiments}\label{sec:numerical_experiments}
    \subsection{Problem setup and data}
        We apply our method to an EUV reflectometry problem using measurement data provided by the ATMOC project \cite{hansen2025atmoceuv, Saadeh_rutheniumsample_2021}. Using our method, we aim to determine the posterior distributions $\pi(x \mid y)$ of the optical constants of the measured specimen, which describe its optical behavior in the EUV range. Aside from the optical constants, we can also provide estimates for the actual thickness of the manufactured specimen, for which only the nominal manufacturing target is known. For our test case, we fix the roughness parameters at the inter-layer interfaces to a small constant value in order to focus on the reconstruction of the main parameters of interest, the optical constants. Evaluation of the forward model \eqref{eq:matrixmethod} with the obtained posterior distributions is then used to validate the obtained results with respect to the physics.

        The MCMC used to solve the inverse problem is based on Bayes' theorem
        
        \begin{align}
            \pi(x \mid y) = \frac{\ell(x \mid y) \pi_0(x)}{\int \ell(x\mid y) \pi_0(x) \, \diff x} \, \mathrm{.}
       \label{eq:bayestheorem}
         \end{align}
         
        Here, $x$ describes all the parameters of interest, while $y$ contains the measurement data resolved at various angles and photon energy levels. The unnormalized posterior $\rho_{\infty}$ described in \autoref{sec:rJKOtt} can be expressed as $$ \rho_{\infty}(x) \propto \ell(x\mid y)\pi_0(x) \, \mathrm{.} $$

        We consider a two-layered specimen with layers of an assumed $ \ce{RuO2} $ contamination layer on top of the (main) $\ce{Ru}$ layer. This contamination layer is added into the model to account for possible impurities and contamination that can occur during both fabrication and handling of the specimen, for example during the actual measurement process. 
        The layer system is situated on top of a  bulk $\ce{Si}$-substrate which is considered to be of infinite thickness for modeling purposes.
        The former assumption is due to the penetration depth of the light used in the measurement being much smaller than the thickness of the bulk material.
        
        Measurements were performed at 11 different wavelengths ($10\mathrm{nm}-20\mathrm{nm}$, spaced at $1\mathrm{nm}$), with each wavelength being resolved at 45 different incident angles ranging from 2\textdegree \, to 88\textdegree. 
        The measurement data provide the measured reflectance values $R$ with $N_{\mathrm{data}}=11\times45 = 495$ data points, as well as an empirical uncertainty budget $\sigma_i$, which is also used during the MCMC procedure. Its derivation is described in detail in \cite{hansen2025atmoceuv,scholze2001highaccuracyeuv}.  
        The nominal thickness of the the investigated $\ce{Ru}$ layer is $30 \mathrm{nm}$.       

        The forward model $f(x)$ described in \eqref{eq:matrixmethod} is used to introduce the measurement data into the MCMC procedure via the logarithmic likelihood function  
        \begin{align}
            \ell_{\text{log}}(x | y) = \sum_{i=1,N_\mathrm{data}}  - \frac{1}{2}\log \left(\frac{1}{{2\pi\sigma_i^2}}\right) - \left(\frac{(y_i-f_i(x))^2}{2\sigma_i^2}\right), \label{eq:definition_of_mcmc_likelyhood}
        \end{align}
        which also incorporates the measurement uncertainties $\sigma_i$. 
        This uncertainty incorporates different sources introduced by the measurement instrumentation, sample variance and modeling errors, as described in \cite{hansen2025atmoceuv}.
        For the proof-of-concept experiments presented in this paper, we have used the measurements at only two wavelengths, namely $10~\mathrm{nm} $ and $ 16~\mathrm{nm}$.
        In such setup the parameter space is already large enough to be challenging for the baseline method, but it is still feasible to perform multiple exploratory experiments.
        Addition of subsequent wavelengths would result in increase of the problem's dimensionality and potentially higher concentration of the posterior.

            To generate the reference sample, we run an ensemble MCMC with $N_{\text{ens}} = 512$ parallel chains. Our setup for the MCMC considers $d=12$ parameters in total, as detailed in ~\autoref{tab:parameters}. All parameters use uniformly distributed priors in their given bounds. Given the small expected thickness of the \ce{RuO2} top layer, the roughness of the top layer interface was included into the forward model as a fixed parameter, with a set value of $0.1 \mathrm{nm}$.

            \begin{table}
                \centering
                \begin{tabular}{c|cc}
                    Parameter & prior bounds & \\ \hline
                    thickness RuO2 layer $[\mathrm{nm}]$ & $\left[ 0.1 , 2 \right]$ & \\
                    thickness Ru layer $[\mathrm{nm}]$ & $\left[ 27, 33 \right]$ & \\
                    roughness RuO2-Ru interface $[\mathrm{nm}] $ &$\left[ 0.5, 1.5 \right]$&\\
                    roughness RU-Si interface $[\mathrm{nm}] $& $\left[ 27, 33 \right]$&\\
                    RuO2 $n(10\mathrm{nm})$ & $\left[ 0.66, 1.26 \right]$&\\
                    RuO2 $k(10\mathrm{nm})$ &$\left[ 0, 0.31 \right]$ &\\      
                    RuO2 $n(16\mathrm{nm})$ &$\left[ 0.59, 1.19 \right]$ &\\
                    RuO2 $k(16\mathrm{nm})$ & $\left[ 0, 0.35 \right]$&\\
                    Ru $n(10\mathrm{nm})$ &$\left[ 0.85, 1.05 \right]$ &\\
                    Ru $k(10\mathrm{nm})$ &$\left[ 0, 0.11 \right]$ &\\      
                    Ru $n(16\mathrm{nm})$ &$\left[ 0.73, 0.93 \right]$ &\\
                    Ru $k(16\mathrm{nm})$ &$\left[ 0, 0.17 \right]$ &\\
                \end{tabular}
                \caption{List of parameters and their priors bounds }
                \label{tab:parameters}
            \end{table}

            The implementation of the MCMC algorithm is provided by the \textit{emcee} Python package \cite{foreman2013emcee}.
            Proposals are generated using the "Stretch Move" described in~\cite{goodman2010ensemble}. For our cases, $  5\cdot 10^5$ samples were generated per chain. The first $\sfrac{1}{6}$ of the samples was discarded as burn-in to not skew the overall results. Convergence was checked by evaluation of the integrated autocorrelation time $\tau_i$ for the affine invariant sampler in \emph{emcee}, as suggested by \cite{goodman2010ensemble}. The computed autocorrelation times for the different parameters were found to be in agreement with the chain lengths of $50 \tau_i$ suggested by \cite{foreman2013emcee}.  The average acceptance rate across all chains is approx. $0.253$, which lies in the suggested range for affine invariant samplers \cite{Roberts01011998}.

         \subsection{Validation of the baseline MCMC}

         To assess the validity of the generated MCMC baseline samples as well as the samples generated via rJKOtt (cf.~\autoref{sec:rJKOtt}), we compare them to the measured data. To this end, we first compute the median values for each parameter of interest from the generated samples and use them to generate a reflectance map using the forward model described in \eqref{eq:matrixmethod}. Comparison is done using the relative $L_2$ norm of the reflectance maps $R_{\mathrm{model}}$ and $R_{\mathrm{measure}}$ , as well as a comparison of the predicted value ranges to the physical expectations. The relative error
         \begin{equation}\label{eq:rel_err}
             \varepsilon_{\mathrm{rel}} = \sqrt{\frac{\sum_{i=1}^{N_{\text{data}}} |R_{\mathrm{measure},i} - R_{\mathrm{model},i}|^2}{\sum_{i=1}^{N_{\text{data}}} |R_{\mathrm{measure},i}|^2}}
         \end{equation}
 of our baseline was in the $10^{-3}$ range when compared to the measurement signal, indicating good agreement to the measurement data.
        
\begin{figure}
            \centering
            \includegraphics[width=1\linewidth]{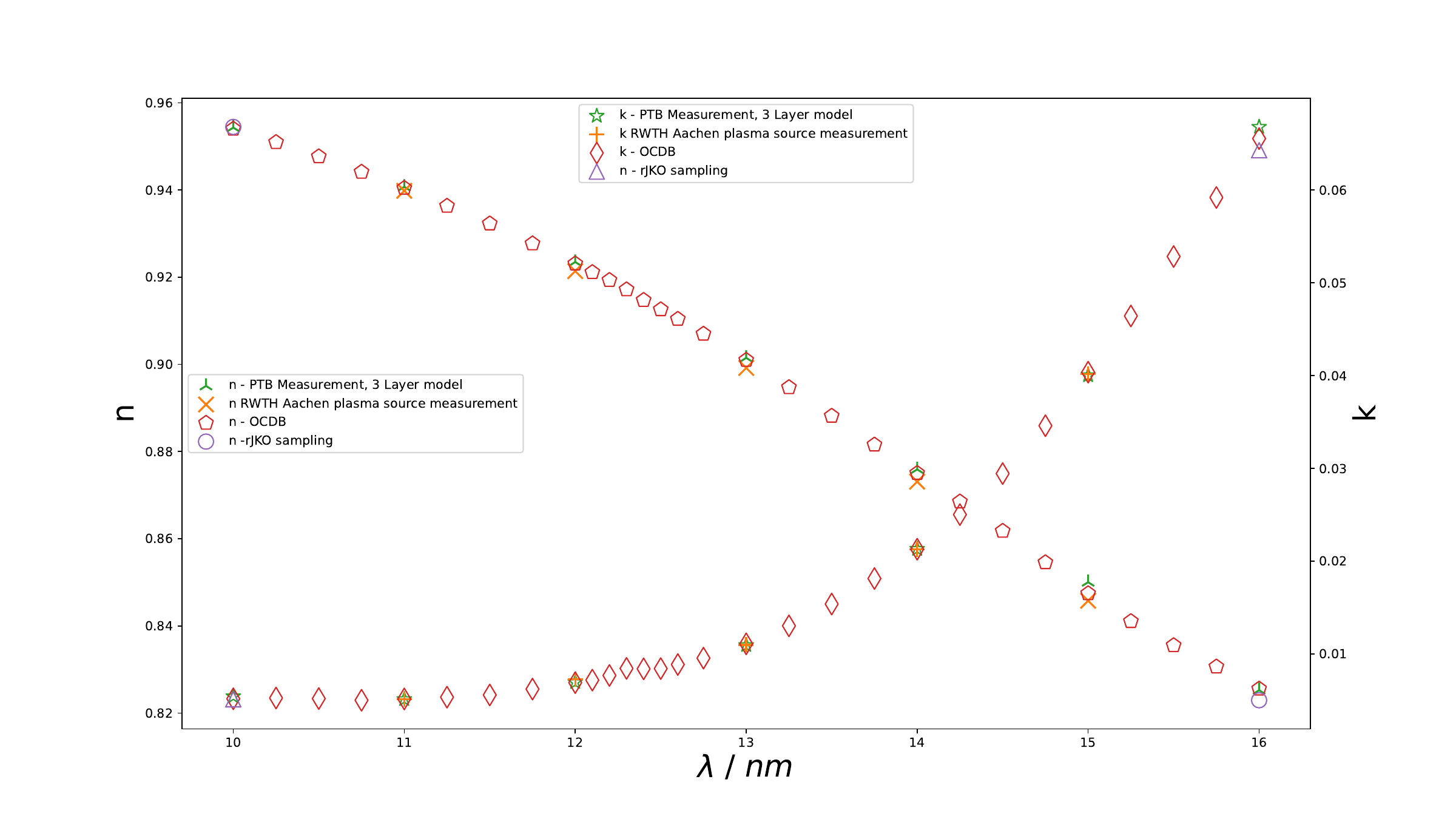}
            \caption{Comparison of the the optical constants for a Ruthenium (RU) layer obtained via different methods. The reference values given by the OCDB \cite{biksup2024ocdb} are depicted in red, the MCMC reconstruction based on a plasma source measurement \cite{hansen2025atmoceuv} is shown in orange. Results of our Baseline MCMC 3 Layer reconstruction are shown in green, exemplary results from rJKOtt are shown in purple. All values depicted are the median values of the sampled distributions.}
            \label{fig:oc_literature_comparison}
        \end{figure}

        For further validation, we compare our baseline MCMC results to reference values given by the OCDB database \cite{biksup2024ocdb} as well as MCMC reconstruction results based on measurements performed the same physical sample specimen investigated in this study. A comparison of the real ($n$) and imaginary parts ($k$) of the optical constants is depicted in \autoref{fig:oc_literature_comparison}. Our results obtained via MCMC lie close to the values given by the OCDB database, as well as to results obtained by MCMC reconstruction using a different measurement technique discussed in \cite{hansen2025atmoceuv}. Note that the OCDB reference data was obtained on a thicker layer of material than the measurements used for this work. It is therefore to be expected that the calculated optical constants do not match perfectly to this reference.

        \subsubsection{Preconditioning and initial approximation}
            We run the Laplace method to estimate mean and covariance of the distribution.
            Finding the MAP estimate in the Laplace method requires solving the minimization problem~\eqref{eq:laplace_MAP}, which we do by means of the COBYQA method~\cite{ragonneau2022model,ragonneau2025cobyqa}.
            The method is a derivative-free version of a trust-region SQP method, relying on quadratic models built with under-determined interpolation of the target functional. 
            We run the method multiple times with different initial points, sampled from the prior bounds, and choose the converged solution with maximal likelihood.
            In our case, finite-difference approximation of~\eqref{eq:laplace_covariance} did not yield a positive-definite estimate with tolerable error in each component.
            Thus we perform iterative importance sampling to improve the estimates  mean and covariance, with the initial covariance estimate being obtained from the prior bounds.
            In total, MAP estimation and importance sampling require approximately $2.5\cdot 10^5$ posterior calls.

            After acquiring mean and covariance estimate for the distribution, the grid is constructed in the transformed coordinates $z$~\eqref{eq:orig_to_pc_transform}.
            The grid is uniform over the interval $[-5, 5]$ (according to the empirical $5\sigma$ rule) with $N_{\text{nodes}} = 30$ nodes in each dimension.
            % [G13 hidden] The grid is presented in~\autoref{fig:tt_init_pc}.
            % [G13 hidden] \begin{figure}
            % [G13 hidden] \centering
            % [G13 hidden] \includegraphics[width=0.5\linewidth]{fig/tt_init_approx_pc_coords.pdf}
            % [G13 hidden] \caption{Initial TT approximation in the transformed $(z)$ coordinates}
            % [G13 hidden] \label{fig:tt_init_pc}
            % [G13 hidden] \end{figure}
            % [G13 hidden] The initial distribution (blue) is a rank-one TT approximation of the standard normal distribution $z \sim \mathcal{N}(0, \Id)$.
            % [G13 hidden] The transformed target distribution (red) is close in the mean and covariances, but exhibits a more complex structure.
            
        \subsection{Learning the posterior distribution}
\begin{figure}
                \centering
                \includegraphics[width=\linewidth]{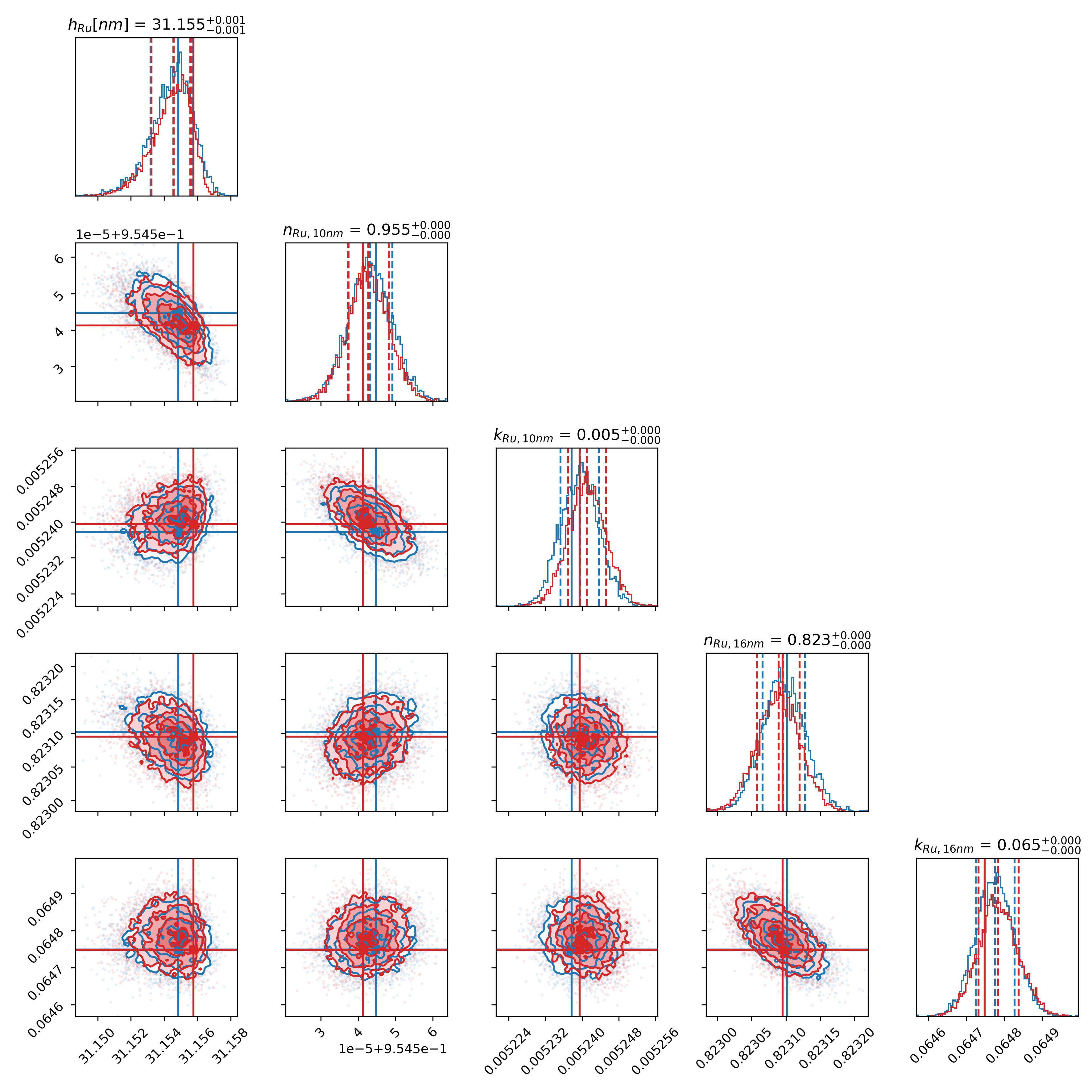}
                \caption{%
                    Posterior samples by the baseline MCMC (blue) and {rJKOtt} (red).
                    One- and two-dimensional marginals of thickness of the main  $\ce{Ru}$ layer are presented, along with maximum likelihood estimates, means and the centered interval containing $68\%$ of the distributions.
                    The shapes of the distributions, as well as their main statistics are well aligned.
                }
                \label{fig:tt_sol_orig_coord}
            \end{figure}

            We perform one step of rJKOtt using regularization factor $\beta_0 = 10^{-1}$ and the effective time in the heat equation $T_0\beta_0 = 10$.
            The maximal rank for the TT approximation is chosen equal to $15$.
            $N_{\text{fp}} = 20$ fixed-point iterations are performed using Anderson mixing method with memory size 1 and relaxation $1.0$.
            The method requires approximately $3.1\cdot 10^6$ posterior evaluations.
            In \autoref{fig:tt_sol_orig_coord}, we show the distribution, given by the TT method (blue) and the MCMC reference (red).
            The TT method, along with the linear preconditioning mapping, is able to adequately capture the non-Gaussian structure of the target and incorporates the prior bounds.

            The total amount of posterior evaluations for both the preconditioning step and the TT solution is approximately $N_{\text{tt}} \approx 3.3 \cdot 10^6$.
            For the ensemble size $N_{\text{ens}} = 512$, this is an equivalent of roughly $6\,500$ iterations, while the reference chain ran for $500\,000$ iterations (total number of posterior calls $N_{\text{MCMC}} \sim 2.6 \cdot 10^8$).
            Thus, the number of posterior calls in our method is significantly less than for the baseline MCMC.
            The computational effort is additionally justified if the TT solution can be re-used in further measurements for a family of similar specimens.
            In the following, we evaluate our method in this setting.

    \subsection{Re-using the solution for data shifts}
   
        In order to leverage the fact that sampling from our rJKOtt approximation is generally faster than sampling with MCMC, we want to explore a use case which makes use of this feature. For this, we consider a process which produces several specimens manufactured to a set of specifications, which are then measured using the reflectometry method outlined in the previous sections. As with every fabrication process, small deviations from the target specifications are generally unavoidable. To model this process, we generate baseline MCMC samples and then shift the obtained distributions in a controlled way to model manufacturing deviations.
  
        Since the generation of samples from the TT approximation is generally faster than a full run of MCMC, we want to demonstrate a use case in which we first compute a TT approximation for the previously described reflectometry problem. Using this as a new starting point, we supply a new, synthetically generated measurement dataset which varies one of the quantities of interest in a controlled manner.

        We generate these shifted synthetic measurements by first solving the inverse problem for the real measurement data for a Ru layer of nominal thickness \qty{30}{nm}, as described in the previous section. The obtained posterior distribution $\pi_{30}(x | y)$ is then shifted in one component, the layer thickness to a new value of $h_{\mathrm{Ru}} = 30.25 \mathrm{nm}$, and the forward model is then applied to the resulting new distribution to obtain a new, synthetic measurement dataset. Additionally, we compute the corresponding potentials $\eta,\ \hat\eta$ using the TT approximation.

        The new dataset is then used to compute the shifted posterior $\rho_{\text{shift}}(x | y_{\text{shift}})$ using rJKOtt, with an accompanying MCMC for validation. Since the input parameter distributions used to generate the synthetic data are fully known, validation against the thus known ground truth is straightforward.

        We repeat the approximation of moments via the Laplace method and importance sampling.
        Then, we restart the fixed-point solution, using the TT approximation for the potentials $\eta,\ \hat\eta$ from the previous case.
        This ``warm-start'' technique improves the performance of the method.
        In our particular case, the fixed-point converges in $n_{\text{fp}} = 4$ iterations and uses $1.2\cdot 10^5$ iterations, which is $\times 10$ less than in the original case.
        \begin{figure}
            \centering
            \includegraphics[width=\linewidth]{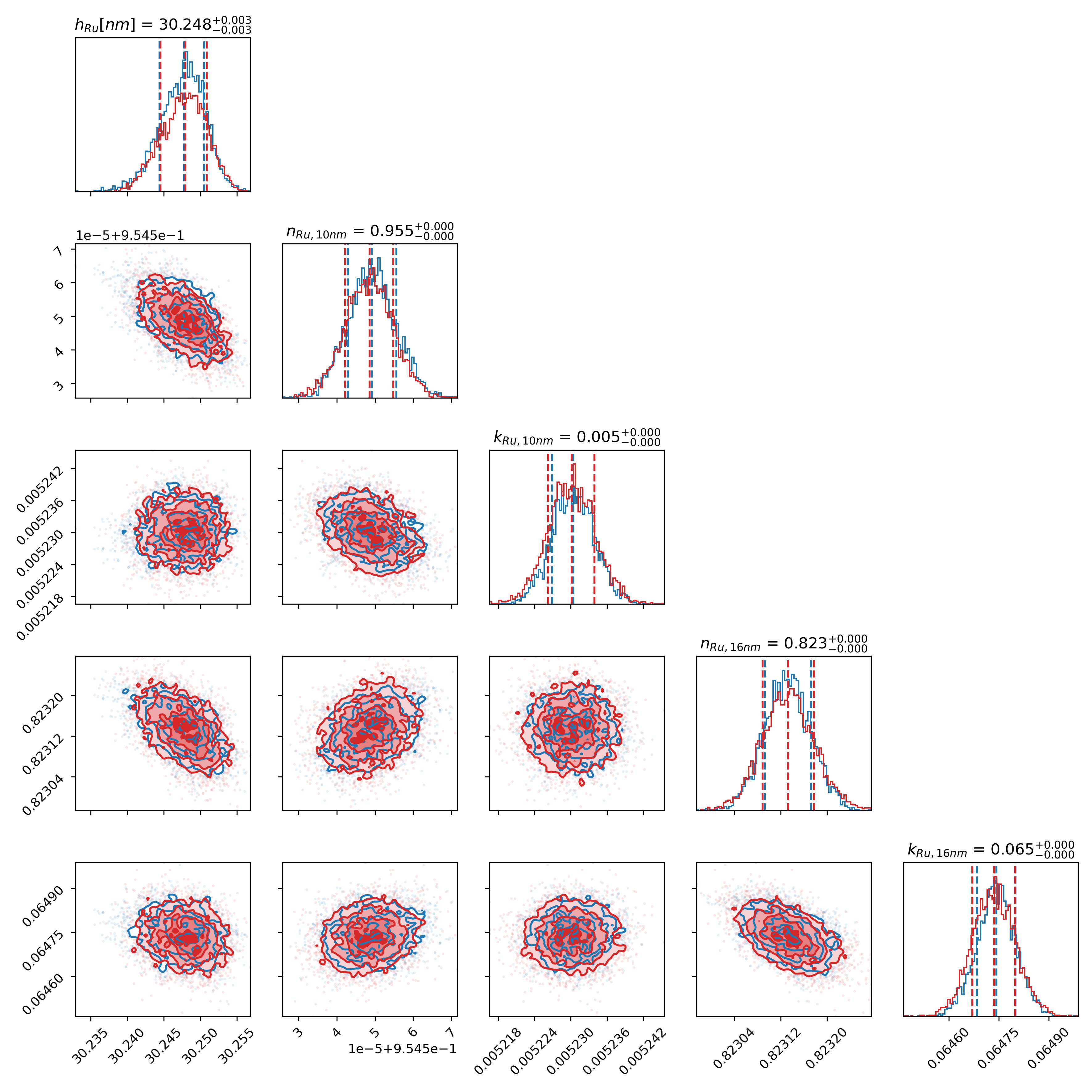}
            \caption{
                    Posterior samples by the baseline MCMC (blue) and {rJKOtt} (red) in the case with shifted thickness.
                    One- and two-dimensional marginals of thickness of the main  $\ce{Ru}$ layer are presented, along with maximum likelihood estimates, means and the centered interval containing $68\%$ of the distributions.
                    The shapes of the distributions, as well as their main statistics are well aligned.
            }
            \label{fig:tt_sol_shift_orig_coord}
        \end{figure}
        The result is presented on~\autoref{fig:tt_sol_shift_orig_coord}. This figure shows the posterior distributions obtained with both MCMC and  the rJKOtt method. For visual clarity, only the distributions for the thickness of the main Ru layer and its optical constants at the two considered wavelengths are depicted, which are also the main quantities of interest in our investigation.
        As can be seen in the figure, the target layer thickness of $h_{\mathrm{Ru}} = 30.25 \mathrm{nm}$ is reflected very well in the obtained posterior distributions, with the median value deviating less than $0.01\%$ for both the rJKOtt as well as the validation MCMC computations.

\section{Conclusion and outlook}

We demonstrated that our method based on rJKOtt in conjunction with a suitable preconditioning can provide samples from the target distribution using less posterior calls than a classic Bayesian MCMC approach without compromising on the quality of the results. 
Even the moderate reduction by a factor of $\sim 10$ described in \autoref{sec:numerical_experiments} can already be considered significant, since this scales directly with the expense of the forward model. Furthermore, we have shown that our method can reuse previous solutions as a starting point to similar problems in order to provide a significant reduction required posterior calls by an order of magnitude.

With our results, we could already demonstrate a reduction of the number of required posterior calls to solve the inverse problem. 
In the future, this methodology could be extended to inverse problems which require more complex forward models, for example PDE solvers based on the finite element methods. 
Compared to the forward model used herein, those generally incur much higher computational costs \cite{truong2025mlgratingfem, gross2012gratingfem}. 
Thus the impact of the saved posterior calls on the overall runtime of the computation will be more pronounced. 

The posterior considered in the present paper, based on measurements from only two wavelengths, already provides a sharp estimate of the specimen's parameters at these wavelengths.
The inclusion of the whole set of measurements remains relevant.
However, it is hindered by the increased dimensionality of the parameter space and concentration of density.
In particular, the estimation of the nonlinear term in~\eqref{eq:fp_terminal_KL} by cross-approximation might be problematic due to larger range of values of $\rho_\infty$, which is harder to represent with finite precision.
The solution of these challenging problems is left to future work.

\paragraph{Acknowledgments}
 
This work has received support from the 20IND04-ATMOC funded by the 
EMPIR programme co-financed by the Participating States and from the European Union’s
Horizon 2020 research and innovation programme.

The authors would like to thank PTB Adlershof group 7.14 for kindly providing the reflectometry measurement data for the $30 \mathrm{nm}$ Ruthenium layer specimen used in this work.

\bibliographystyle{alpha}
\bibliography{references}
\end{document}